\documentclass[10pt]{article}
\usepackage{amssymb,amsbsy,amsmath,amsfonts,amssymb,amscd}
\usepackage{latexsym}

\begin{document}
\def\ddd{\displaystyle}
\def\C{{\mathbb C}}
\def\N{{\mathbb N}}
\def\Z{{\mathbb Z}}
\def\R{{\mathbb R}}
\def\T{{\mathbb T}}
\def\k{{\bf k}}
\def\Ext{{\rm Ext}}
\def\ma{{}_AM}
\def\ss{\langle X \rangle}
\def\sv{\k\langle X \rangle}
\def\bb{\overline b}
\def\epsilon{\varepsilon}
\def\phi{\varphi}
\def\pphi{\widetilde\varphi}
\def\kappa{\varkappa}
\def\wz{\thinspace}
\def\proof{P\wz r\wz o\wz o\wz f.\hskip 6pt}
\def\quest#1{\hskip5pt {\scshape  Problem} {\rm #1}.\hskip 6pt}
\def\leq{\leqslant}
\def\geq{\geqslant}
\def\pd#1#2{\frac{\partial#1}{\partial#2}}
\def\limsup{\mathop{\overline{\hbox{\rm lim}\,}}}
\def\ug#1#2{\left\langle#1,#2\right\rangle}
\def\ii{{\mathrel{\Longrightarrow}}}
\font\SY=msam10
\def\square{\hbox{\SY\char03}}

\def\C{{\mathbb C}}
\def\N{{\mathbb N}}
\def\Z{{\mathbb Z}}
\def\R{{\mathbb R}}
\def\PP{\cal P}
\def\p{\rho}
\def\phi{\varphi}
\def\ee{\epsilon}
\def\ll{\lambda}
\def\l{\lambda}
\def\a{\alpha}
\def\as{\sigma}
\def\b{\beta}
\def\D{\Delta}
\def\g{\gamma}
\def\rk{\text{\rm rk}\,}
\def\dim{\text{\rm dim}\,}
\def\ker{\text{\rm ker}\,}
\def\square{\vrule height6pt width6pt depth 0pt}
\def\epsilon{\varepsilon}
\def\phi{\varphi}
\def\kappa{\varkappa}
\def\proof{P\wz r\wz o\wz o\wz f.\hskip 6pt}
\def\quest#1{\hskip5pt {\scshape  Problem} {\rm #1}.\hskip 6pt}
\def\leq{\leqslant}
\def\geq{\geqslant}
\def\pd#1#2{\frac{\partial#1}{\partial#2}}
\def\limsup{\mathop{\overline{\hbox{\rm lim}\,}}}
\def\ug#1#2{\left\langle#1,#2\right\rangle}
\def\kk#1#2{{\k\langle#1,#2\rangle}}
\def\sv{\bf{ k} \langle X \rangle}
\def\k{{\bf k }}
\def\lxr{\langle X \rangle}
\def\defin#1{\smallskip\noindent
{\scshape  Definition} {\bf #1}{\bf .}\hskip 8pt\sl}

\newtheorem{lemma}{Lemma}[section]
\newtheorem{theorem}[lemma]{Theorem}
\newtheorem{corollary}[lemma]{Corollary}
\newtheorem{proposition}[lemma]{Proposition}
\newtheorem{definition}[lemma]{Definition}

\title{Graphs of relations and Hilbert series }

\author{Peter Cameron \&
 Natalia Iyudu}

\date{}

\maketitle

\small

\centerline{Queen Mary, University of London}

\smallskip

\centerline{ {\bf e-mail:} \ n.iyudu@qmul.ac.uk}

\bigskip

\small

{\bf Abstact} We are discussing certain combinatorial and counting
problems related to quadratic algebras. First we give examples which
confirm the Anick conjecture on the minimal Hilbert series for
algebras given by $n$ generators and $\frac {n(n-1)}{2}$ relations
for $n \leq 7$. Then we investigate combinatorial structure of
colored graph associated to relations of RIT algebra. Precise
descriptions of graphs (maps) corresponding to algebras with maximal
Hilbert series are given in certain cases. As a consequence it turns
out, for example, that RIT algebra may have a maximal Hilbert series
only if components of the graph associated to each color are
pairwise 2-isomorphic.

\bigskip


\normalsize\hfill\break \rm

\vskip.5cm

\hrule
 \tableofcontents

\vskip.5cm

\bigskip

\section{Introduction}

Let ${\cal A}(n,r)$ be the class of all graded quadratic algebras on
$n$ generators and $r$ relations:
$$
A= k\langle x_1,\dots,x_n\rangle / {\rm id}\{p_i:i=1,\dots,r\},
$$
where $p_i=\sum\limits_{k,j=1}^{n} \alpha_i^{k,j} x_k x_j,$
$\,\,\alpha_i^{k,j} \in k$.

We deal with an arbitrary field $k$ of char $0$. Only on the way
(section \ref{gp}) we restrict ourself with $\mathbb C$ for a while
(to get more general statement), but it will not influence further
results.

These algebras are endowed with a natural filtration
$A=\bigcup\limits_{m=0}^\infty U_m$, where $U_m$ is the linear span
of monomials on $a_{i}$ of degree not exceeding $m$, $a_i$ are the
images of the variables $x_i$ under the canonical map from $k\langle
x_1,\dots,x_n\rangle$ to $A$ and the degree of $a_{i_1}\dots
a_{i_d}$ equals $d$. Since the generating polynomials are
homogeneous, the algebra $A\in \cal A$ also possesses   a canonical
grading $A=\bigoplus\limits_{i=0}^\infty A_i$, where $A_i$ is the
linear span of monomials of degree exactly $i$. This grading has a
finiteness property: ${\rm dim}_k A_i < \infty $ for any $i$, since
algebra is finitely generated. This allow to associate to the series
of dimensions various generating functions. The one which reflects
most straightforward properties of the algebra will be considered
here.
\begin{definition}\label{H} \ The Hilbert series of a graded algebra
$A=\bigoplus\limits_{i=0}^\infty A_i$ is the generating function of
the series of dimensions of graded components $d_i={\rm dim}_k
\,A_i$ of the following shape: $ H_A(t)=\sum_{i=0}^\infty d_it^i. $
\end{definition}

We are going to confirm  Anick's conjecture \cite{An} saying that a
lower bound for the Hilbert series of an algebra with
$\frac{n(n-1)}{2}$ quadratic relations given by the series
$\left|\left(1-nt+\frac{n(n-1)}2t^2\right)^{-1}\right|$ is attained,
for the small number of variables $n \leq 7$. Here the sign of
modulus stands for the series where $n$th coefficient equals to the
$n$th coefficient of initial series if this is positive together
with all previous coefficients and is zero otherwise.

After notices on  minimal and generic series for quadratic algebras
we turn to the main subject of our investigation. We consider
subclass ${\cal R}(m,n) \subset {\cal A}(g,\frac{g(g-1)}{2})$, where
$g=n+m$, called RIT algebras (it was introduced and studied in
papers \cite{a,ia,wia}). Class ${\cal R}(m,n)$ defined as consisting
of algebras with presentation of the form

$$R=k \langle x_1,...,x_m,y_1,...,y_n \rangle / F,$$

where

\begin{equation}
F=id \left\{
\begin{array}{l}
{[x_i,x_j]=0, }
\\
{[y_i,y_j]=0, }
\\
{[x_i,y_j]= y_{f(i,j)} y_{j},}
\end{array}
\right.\label{fff}
\end{equation}
\noindent and $f$ is a map $ f:M \times N \longrightarrow N, \,
M=\{1,...,m\}, \, N=\{1,...,n\}$.

To any such algebra we associate $m$-colored graph with $n$ vertices
in such a way that subgraph $\Gamma_i$ of color $i$ reflects the map
$\sigma_i: N \to N$ defined by $\sigma_i(n)=f(i,n)$.

We formulate conditions (Theorem \ref{GrSig}) on the above maps
which mean that defining relations of algebra form  a Gr\"obner
basis, or equivalently that the algebra has a lexicographically
maximal Hilbert series. Then we attack more subtle question on how
to describe precisely the combinatorial structure of those maps.
This is done explicitly in the Theorem \ref{s1s2} for a pair of
maps. As a consequence  interesting necessary condition is obtained
for an algebra to have a maximal Hilbert series. Corollary
\ref{cs1s2} says that Hilbert series of algebra could be maximal
only if graphs of all maps $\sigma_i, i=\bar{1,m}$ are 2-isomorphic.
We call graphs {\it 2-isomorphic} if they become isomorphic after
gluing pairs of vertices in common cycles of length two.

Described combinatorial conditions (Theorem \ref{s1s2}) also imply
that algebra is Koszul and obeys a commutator generalization of
Yang-Baxter equation. Another consequences for RIT algebras from
being presented by a quadratic Gr\"obner basis are that in this case
they  are Auslander regular and Cohen-Macaulay. Hence we get
combinatorial conditions on graphs sufficient also for obeying these
properties.

In the section   \ref{tc} we present the complete list of Hilbert
series and corresponding non isomorphic colored graphs for RIT
algebras of rank up to 4.

\section{The Anick conjecture for $n \leq 7$}

\subsection{Series in general position}\label{gp}

We remind here the proof of minimality of general series because it
is essential for the next section. Version we present deals with the
notion of general position in the Lebesgue sense, so we suppose for
this section that the field is $k=\mathbb C$.
Essentially the knowledge on this matter is due to Anick \cite{An}
and explanations in simplified  form one can find in the survey of
Ufnarovskii \cite{Ufn}. In \cite{PP} one can also find a remark on
the minimality of $n$th component of the series in general position
in the  Zarisski sense. But as it is pointed out in \cite{Ufn}
(remark before the theorem 3 in I.4.2), since the infinite union of
proper affine varieties may not be contained in a proper affine
variety, one can not state  that the minimal series (in case it is
infinite) is in general position in the Zarisskii sense. There could
be several ways to avoid this problem, for example in spirit of
theorem 3, I.4.5. in \cite{Ufn}. But over $\mathbb C$ most natural
and easy way is to use the topology defined by the Lebesgue measure,
so we present this version here.

 Let we define now more precisely what is meant by an algebra in
general position in the Zarisskii and in the Lebesgue sense.
Algebras from ${\cal A}(n,r)$ are naturally labeled by  the points
of $k^{rn^2}$, corresponding to the coefficients $\alpha_i^{k,j}$ of
the relations. Given a property $P$ of quadratic algebras, we say
that $P$ {\it is satisfied for $A\in {\cal A}(n,r)$ in general
position in the Zarisskii sense} if the set of the coefficient
vectors corresponding to those $A\in {\cal A}(n,r)$, which obey the
property $P$, is a non-empty Zarisskii-open subset of $k^{rn^2}$. We
also say that $P$ {\it is satisfied for $A\in {\cal A}(n,r)$ in
general position in the Lebesgue sense} if the set of the
coefficient vectors corresponding to those $A\in {\cal A}(n,r)$,
which do not obey the property $P$ has $rn^2$-dimensional Lebesgue
measure zero. Since the set of zeros of any non-zero polynomial has
the Lebesgue measure zero, we see that as far as arbitrary property
$P$ satisfied for $A\in {\cal A}(n,r)$ in general position in the
Zarisskii sense it is also satisfied for $A\in {\cal A}(n,r)$ in
general position in the Lebesgue sense. Define the minimal series in
the class ${\cal A}(n,r)$ componentwise:
$$
H^{n,r}_{\rm min}(t)=\sum_{i=0}^\infty b_it^i,
$$
where $b_i=\min\limits_{A\in {\cal A}(n,r)}{\rm dim}\,A_i$, $A_i$
being the $i$th homogeneous component in the grading of $A$. It is
not clear {\it a priori}, whether there exists an algebra $A\in
{\cal A}(n,r)$ whose Hilbert series coincides with $H^{n,r}_{\rm
min}$. This follows however from the statement below.

\begin{proposition}\label{t1} For $A\in {\cal A}(n,r)$ in general
position in the Lebesgue sense, the equality $H_A=H^{n,r}_{\rm min}$
is satisfied.
\end{proposition}

\begin{proof}
 Denote the ideal generated
by $\{p_i:1\leq i\leq r\}$ by $I$ and its $d$th homogeneous
component by $I_d$. Obviously
$$
I_d={\rm span}_k\,\{up_iv:u,v\in \langle x_1,\dots,x_n\rangle,\,\,
{\rm deg}\,u+{\rm deg}\,v=d-2\}.
$$
Here $\langle x_1,\dots,x_n\rangle$ stands for the free semigroup
generated by $\{x_1,\dots,x_n \}$. Let $w_1,\dots,w_m$ be all
monomials of degree $d$ in the free algebra $k\langle
x_1,\dots,x_n\rangle$. Since it is a linear basis in the $d$th
homogeneous component of $k\langle x_1,\dots,x_n\rangle$, we can
uniquely express the above polynomials $up_iv$ as a linear
combinations of $w_l$:
$$
up_iv=\sum_{l=1}^m \lambda_{u,v,i}^l w_l.
$$
The dimension of $I_d$ is exactly the rank of the rectangular matrix
$\Lambda=\{\lambda_{u,v,i}^l\}$, whose rows of length $m$ are
labeled by the triples $(u,v,i)$, where $i=\bar{1,r}$ and $u,v$ are
monomials in $x_1,\dots,x_n$ satisfying ${\rm deg}\,u+{\rm
deg}\,v=d-2$.

Obviously $\lambda_{u,v,i}^l$ are linear functions of the
coefficients $\alpha_i^{k,l}$ of the polynomials $p_i$. Condition
that the dimension of $A_d$ is minimal is equivalent to the
condition that ${\rm dim}\,I_d={\rm rk}\,\Lambda$ is maximal. Denote
the maximal rank of $\Lambda$ by $D$. Thus, the dimension of $A_d$
is minimal if and only if there is a non-zero minor of the matrix
$\Lambda$ of the size $D$. The family of the minors of $\Lambda$ of
the size $D$ is a finite family of polynomials ${\cal P}_l$ on the
coefficients $\alpha_i^{k,l}$ and some of these polynomials are
non-zero. This means that the set of $A\in{\cal A}(n,r)$ with
minimal ${\dim}\,A_d$ corresponds to the complement of the union of
the sets of zeros of finitely many non-zero polynomials. Any such
set is a non-empty Zarisskii open set and its complement has zero
Lebesgue measure. The set of algebras $A\in{\cal A}(n,r)$ satisfying
$H_A=H^{n,r}_{\rm min}$ is then a countable intersection of
non-empty Zarisskii open sets and therefore its complement has zero
Lebesgue measure as a countable union of sets with the Lebesgue
measure zero. This completes the proof of  the proposition.
\end{proof}\square

{\bf Remark } \ Let us mention that in  case when the minimal series
$H^{n,r}_{\rm min}$ is finite, the countable union from the proof of
the Proposition~\ref{t1} is in fact finite and the equality
$H_A=H^{n,r}_{\rm min}$ is satisfied for $A\in{\cal A}(n,r)$ in
general position in the Zarisskii sense as well and over an
arbitrary field.

\subsection{The Anick conjecture holds for $n \leq 7$}

Now we are back to arbitrary basic field $k$ of char $0$.  We
consider the question whether the minimal series is finite for the
case $r=n(n-1)/2$. It was raised in the paper of D.~Anick \cite{An},
where a lower bound for the Hilbert series for algebras from ${\cal
A}(n,\frac{n(n-1)}{2})$ was discovered. It was established that
$$
H^{n,n(n-1)/2}_{\rm
min}\geq\left|\left(1-nt+\frac{n(n-1)}2t^2\right)^{-1}\right|,
$$
where $\geq$ is a componentwise inequality, i.e. it holds if each
coefficient of the first series is greater or equal then  the
corresponding coefficient of the second one and $|f(t)|$ stands for
the positive part of the series $f\in k[[t]]$. More precisely, if
$f(t)=a_0 + a_1 t +  a_2 t^2 \dots$, then $|f(t)|=b_0 + b_1 t +  b_2
t^2 \dots$, where $b_m=a_m$ for $m \in \{i\,|\, a_j \geq 0
\,\,\forall j \leq i\}$ and $b_m=0$ otherwise. There was a question
raised whether this lower bound is attained.

Since we know from the theorem \ref{t1} that the algebras of ${\cal
A}(n,r)$ in general position have minimal Hilbert series, to prove
that this estimate attained it would be enough to be able to write
down generic coefficients of the relations and calculate the Hilbert
series.

{\bf Example 1.} \ The algebra $A$ over the field $k=\mathbb Z_{17}$
given by the relations
$$
A=k\langle a,b,c\rangle\Big/\left\{\begin{array}{l}
ac+2ba+9b^2+3ca+9cb+8c^2,\\
3ab+5ac+7ba+b^2+8bc+4ca+cb+2c^2,\\
10a^2+2ab+11ac+2ba+8b^2+4bc+9ca+7cb+5c^2
\end{array}\right.
$$
has the Hilbert series $H_A=1+3t+6t^2+9t^3+9t^4=|(1-3t+3t^2)^{-1}|$.
\footnote {The computations were done using bunch of programs GRAAL
(Graded Algebras) written in Uljanovsk by A.Kondratyev under the
guide of A.Verevkin.}

By this method we are able to confirm Anick's conjecture for small
number of generators.

\begin{proposition}\label{t2}
The lower bound for the Hilbert series of an algebra $A \in { \cal
A} (n,\frac{n(n-1)}{2})$ over a field $k$ of char $0$ given by
$\left|\left(1-nt+\frac{n(n-1)}2t^2\right)^{-1}\right|$ is achieved
for $n \leq 7$.
\end{proposition}

\proof In Example 1 we have been calculating over the field $k$ of
characteristic $p=17$. Since the series $|(1-3t+3t^2)^{-1}|$ which
is known to be the lower bound coincide with the result of our
calculations, we actually have shown that for any term of the
series, rank of matrix $\Lambda=\{\lambda_{u,v,i}^l\}$ formed as
above in the proof of the proposition \ref{t1}, with
$\lambda_{u,v,i}^l \in \mathbb Z_p$ is maximal. We now  can see that
rank of the same matrix considered over $k$ is also maximal. Indeed,
the reduction from  $\mathbb Z_p$ to $\mathbb Z$ we have because
${\rm rk} M(\mathbb Z_p) \leq {\rm rk} M(\mathbb Z)$. Then since
${\rm char} k =0$, we have $\mathbb Z$ embedded into $k$ and the
same matrix has maximal rank over $k$. So, if over the field
$\mathbb Z_p$ for some $p$, the rank is maximal, then it is maximal
over $k$.

Similarly, we have got an examples of algebras with the Hilbert
series $1+nt+\frac{n(n+1)}{2} t^2 + n^2 t^3 +
n^2(n-\frac{(n^2-1)}{2}) t^4$ for $n=4$ and $1+nt+\frac{n(n+1)}{2}
t^2 + n^2 t^3 $ for n=5,6,7, which are coincide  with the series $|
(1+nt+\frac{n(n-1)}{2} t^2)^{-1}|$ for these values of $n$. \square

\section{RIT algebras and maps of the finite set}

\subsection{The class of RIT algebras}

Here we consider a subclass $\cal R$ of the above class of quadratic
algebras ${\cal A}(n,\frac{n(n-1)}{2})$. The class $\cal R$  of RIT
(relativistic internal time) algebras consists of homogeneous
finitely generated quadratic algebras given by relations of  type
(\ref{fff}).

It turned out that if relations form a Gr\"obner basis, the algebras
from $\cal R$ are so called "geometric rings", more precisely they
are Auslander regular,  Cohen--Macaulay, ${\rm gl dim}\,R={\rm GK
dim}\,R=n$ (number of generators) for them. We have proved this in
\cite{wi} using combinatorial techniques related to the notion of an
$I$-type algebra introduced by J.Tate and M.Van den Bergh in
\cite{TB}. These arguments appeared due to inspiring question  of
M.Van den Bergh on whether RIT algebras obey these properties. It
becomes clear later on that Auslander regularity {\it et.al.} could
also be proved without employing the I-type property, in more
general context, using arguments involving associated graded
structures with respect to appropriate grading (see \cite{Li},
\cite{LevAR}).

The origin of the class  of RIT algebras could be described in such
a way.
 The Lie algebra RIT was introduced
in \cite{a} as a modification of the Poincare algebra ${\cal
P}_4={\cal L}_4+{\cal U}$. Here ${\cal L}_4=O(3,1)$ is a Lorenz
algebra but the space $\cal U$ changed by addition of a new variable
$T$ (related to the relativistic internal time) to the set of
initial variables. The corresponding commutation relations
containing $T$ were derived. Taking an enveloping algebra of ${\cal
P}_4$ and considering the associated graded algebra we obtain the
associative RIT algebra which gives rise to the class under
consideration.

Let we mention that the simplest algebra from this class
$R=R_{1,1}={ k}\langle x,y\rangle\big/ (xy-yx-y^2)$ is one of the
two Auslander regular algebras of global dimension 2,  the second
one is the usual quantum plane ${ k}\langle x,y\rangle\big/
(xy-qyx)$ (this follows from the Artin, Shelter classification
\cite{AS}). We have been studying finite dimensional representations
of it in \cite{prI}.

\subsection{Condition on maps and Gr\"obner basis }

The presentation (\ref{fff}) of the algebra gives us a set of maps
$\sigma_i:N\to N$ defined as
$$
\sigma_i(j)=f(i,j)\ \ \forall j\in N=\{1,\dots,n\},\ i\in
M=\{1,\dots,m\}.
$$

We are interested to relate the properties of algebras to the
properties of these maps.

In particular, we will clarify combinatorial conditions on the
associated colored graph (=set of maps)  which mean that relations
form a Gr\"obner basis. This gives at the same time a condition
equivalent to the maximality of the Hilbert series. Here we mean a
lexicographical order on the series.

To speak about the Gr\"obner basis we have to fix an ordering on the
set of variables, let $x_i>y_j$ for any $i,j$ and $x_i>x_j$ for
$i>j$, $y_i>y_j$ for $i>j$. On the monomials of variables $x_i, y_j$
the order is supposed to be degree-lexicographical.

Then we rewrite relations (\ref{fff}) in the form

\begin{equation}
F=\left\{
\begin{array}{l}
{[x_i,x_j]=0, \, \forall i>j}
\\
{[y_i,y_j]=0, \, \forall i>j}
\\
{[x_i,y_j]= |y_{f(i,j)} y_{j}|, \, \forall i,j}
\end{array}
\right.\label{fff>}
\end{equation}

Here $|y_{f(i,j)} y_{j}|$ stands for the normal form of this
monomial, i.e.

\begin{equation}
|y_{f(i,j)} y_{j}|=\left\{
\begin{array}{l}
{y_{f(i,j)} y_{j} \,\  {\rm if} \ f(i,j) \leq j}
\\
{ y_{j}y_{f(i,j)} \,\ {\rm if}\ f(i,j) > j}
\end{array}
\right.\label{yy>}
\end{equation}

\begin{theorem}\label{GrSig}
The relations of the type (\ref{fff>}) form a reduced Gr\"obner
basis if and only if the function $f(i,j)$ defines the set of
actions $\sigma_i$ with the following property. For any pair of maps
$\sigma_i,\sigma_k, \, k>i$ one of the two conditions is satisfied
in each point $j\in N$:
\\
 either
(1). $\sigma_k(j)=\sigma_i(j)$ and
$\sigma_i\sigma_k(j)=\sigma_k\sigma_i(j)$
\\
or  (2). $ \sigma_k(j) = \sigma_k\sigma_i(j)$ and $ \sigma_i(j) =
\sigma_i\sigma_k(j)$.
\end{theorem}

Let we mention that the second condition implies kind of strong
version of braid type relations on  $\sigma_i$:
$\sigma_k=\sigma_k\sigma_i\sigma_k$ and
   $\sigma_i=\sigma_i\sigma_k\sigma_i$
for $k>i$.

\proof The proof is a direct application of Gr\"obner bases
technique due to Buchberger \cite{Buch} and Bergman \cite{Bergm}.

To find out that relations form a Gr\"obner basis we have to check
that all ambiguities are solvable. Possible ambiguities
in our case are of four types:\\
 1. $x_i x_j x_k, \, i>j>k$, \\
 2. $y_i y_j y_k, \,i>j>k$, \\
 3. $x_i y_j y_l, \, \forall i, j>l$, \\
 4. $x_l x_i y_j, \, \forall j, l>i$.

The first two types are trivially solvable. Ambiguities of the type
three are also solvable:

\def\double#1#2{\hbox{$\vcenter{\offinterlineskip\halign
{##\hfil\cr\hfil #1\cr \vrule height 4pt depth0pt width0pt\cr
#2\cr}}$}}
\def\triple#1#2{\underline{#1}\!\!\!\!\underline{\vrule height
0pt depth5.8pt width0pt \,\,\,\,#2}}
\def\Underline#1{\underline{#1\!}\,}

$$
\begin{array}{rcl}
y_j\underline{x_iy_l}+y_{(i,j)}y_jy_l&\longleftarrow\
\triple{x_iy_j}{y_l}\ \longrightarrow&
\Underline{x_iy_l}y_j\\ \downarrow\qquad\qquad&&\quad\downarrow\\
y_jy_lx_i+y_jy_{(i,l)}y_l&&y_l\Underline{x_iy_j}+y_{(i,l)}y_ly_j\\
&&\quad\downarrow\\
&&y_ly_jx_i+y_ly_{(i,j)}y_j
\end{array}
$$

Ambiguities of the type four are solvable if and only if the
following two-element (non-ordered) sets are coincide:

$$ \{ f(i,j); f(l,f(i,j)) \} = \{ f(i,f(l,j)); f(l,j) \},$$
\noindent for any $j$ and $l>i$.

Indeed:
$$
\begin{array}{rcl}
x_i\underline{x_ly_j}&\longleftarrow\ \triple{x_lx_i}{y_j}\
\longrightarrow&
\Underline{x_ly_j}x_i+\Underline{x_l|y_{(i,j)}}y_j|\\
\downarrow\qquad\qquad&&\quad\downarrow\\
\Underline{x_iy_j}x_l+\Underline{x_i|y_{(l,j)}}y_j|&&\!\!\!\!\!\!\!
\double{$y_jx_lx_i+ y_{(l,j)}y_jx_i+y_{(i,j)}\underline{x_ly_j}+$}
{$+y_{(l,(i,j))}y_{(i,j)}y_j$}\\
\downarrow\qquad\qquad&&\quad\downarrow\\
\double{$y_jx_ix_l+y_{(i,j)}y_jx_l+y_{(l,j)}\underline{x_iy_j}+$}
{$+y_{(i,(l,j))}y_{(l,j)}y_j$}&&\!\!\!\!\!\!\!
y_{(i,j)}y_jx_l+y_{(i,j)}y_{(l,j)}y_j\\
\downarrow\qquad\qquad&&\\
y_{(l,j)}y_jx_i+y_{(l,j)}y_{(i,j)}y_j&&
\end{array}
$$

In some places above we write for example $|y_{(i,j)}y_j|$ in stead
of $y_{(i,j)}y_j$. These are those places where order on $y_i$
essential for the future reductions (namely we have some $x_j$
before $y_i$). We actually had to check all possibilities for the
pairs $|y_{(l,j)}y_j|$, $|y_{(i,j)}y_j|$ appearing at the above
sequences of reductions and all of them via different cancelations
gave the same result.

The coincidence of above mentioned sets means in the language of
maps $\sigma_i$ that
$$\{ \sigma_i(j); \sigma_l\sigma_i (j) \} = \{ \sigma_i\sigma_l(j); \sigma_l(j) \},$$
\noindent for any $j$ and $l>i$.

These sets are coincide if and only if for any fixed $l>i$ in each
point $j$ we have either $ \sigma_i(j) = \sigma_l(j)$ and
$\sigma_l\sigma_i (j)= \sigma_i\sigma_l(j)$ or
$\sigma_i(j)=\sigma_i\sigma_l(j)$ and $\sigma_l\sigma_i
(j)=\sigma_l(j)$. By this we are done. \square

It is of course very natural and important question, when given
presentation of an algebra  form a Gr\"obner basis. In RIT case
these conditions take a specific shape of description of defining
maps $\sigma_i$, obtained above. Conditions for that were formulated
also for example for the class of G-algebras in \cite{LevG} under
the name of non--degeneracy condition. Now we turn to more difficult
matter of clarifying a precise combinatorial structure of  maps
obeying conditions of the theorem \ref{GrSig}. As a first step we
consider few particular cases, which we will  use later on to prove
the general fact.

\subsection{Representations of the semigroup $\langle x_i | x_i=x_i x_j
\rangle $}

Here we consider the case when all elements of $N$ obey conditions
(2) from the theorem \ref{GrSig}. That is we have $ \sigma_k(j) =
\sigma_k\sigma_i(j)$ and $ \sigma_i(j) = \sigma_i\sigma_k(j)$ for
any $k>i, j \in N$. This means that  $\sigma_i$s form  a
representation by actions on the finite set of the semigroup $\Omega
= \langle x_i | x_i=x_i x_j, 1\leq i \neq j \leq m \rangle $. From
these relations it follows that all $x_i$ are idempotents. Reduction
of the first subword $x_ix_j$ in $x_ix_jx_i$ gives
$x_ix_jx_i=x_ix_i$, of the second one: $x_ix_jx_i=x_ix_j$, but then
$x_ix_j=x_i$. Hence we could also write relations just like $\Omega
= \langle x_i | x_i=x_i x_j, i,j=\bar{1,m} \rangle$, without the
condition $i \neq j$. Note, that this semigroup consists in fact of
$m+1$ elements: any word in this semigroup is equal to its first
letter.

What is the structure of maps which form representations then?

\begin{theorem}\label{reps}
Any representation of the semigroup $\Omega=\langle
x_i|x_ix_j=x_i\rangle$ has the following structure. The set of
representation $N$ is decomposed into a disjoint union of subsets.
In each of them there are $m$ fixed points (not necessarily
different), such that the maps $\as_k$, $k=\bar{1,m}$ send the
entire subset to the $k$th of these points.
\end{theorem}

\proof
 Let $\{\as_k\}_{k=1}^m$ be a representation of the semigroup
$\Omega$ on the set $N$. That is, $\as_j\as_k=\as_j$ for any $1\leq
j,k\leq m$. First note that since $\as_k$ are idempotents, they are
identity on their images: $\as_k(m)=m$ for each $m\in R_k={\rm
Im}\,\as_k$. Define the equivalence relation on $N$ corresponding to
$\as_1$: $m_1\sim_{\as_1} m_2$ if $\as_1(m_1)=\as_1(m_2)$. Then $N$
splits into the union of equivalence classes ${\cal
O}_{s_1},\dots,{\cal O}_{s_k}$, where each class contains a unique
element $s_j$ from $R_j={\rm Im}\,\as_j$, so we can enumerate these
classes by these elements. Consider the restriction of the maps
$\as_k$ to an arbitrary class ${\cal O}_{s_j}$. Since
$\as_j\as_k=\as_j$, each $\as_k$ leaves the set ${\cal O}_{s_j}$
invariant. Indeed, $\as_j(\as_k(r))=\as_j(r)=s_j$ for each $r\in
{\cal O}_{s_j}$ and therefore $\as_k(r)\in {\cal O}_{s_j}$.
Moreover, since $\as_k\as_j=\as_k$, the set $\as_k({\cal O}_{s_j})$
consists of one element $\as_k(s_j)$. Indeed,
$\as_k(r)=\as_k(\as_j(r))=\as_k(s_j)$ for each $r\in {\cal
O}_{s_j}$. Hence the structure of these maps is the following: the
set $N$ is decomposed into a disjoint union of subsets, in each of
which $m$ points are chosen and the maps $\as_k$ map the entire
subset to the $k$th of these points. Some of these points could
coincide. \square

Obviously, the other way around, if one take any set of maps
$\{\as_i\}_{i=\bar{1,m}}$ with the described structure, then they
form a representation of the semigroup $\Omega$, that is satisfy the
relations $\as_k\as_i=\as_k, \forall k,i=\bar{1,m}$.

Thus there exists 1-1 correspondence between representations of
$\Omega$ on a finite set and maps described in the theorem
\ref{reps}.

Let we mention that the same is true for representations on an
infinite set, our arguments work there without any change.

It is natural to ask when there exists a faithful representation.

\begin{corollary}\label{n>=m}
For any $n \geq m$ there exists a faithful representation of
$\Omega$ on the set of size $n$.
\end{corollary}

\proof The image of the semigroup ${\Omega}$ in the set of maps
consists just of the maps $\as_1,...,\as_m$, which are images of the
generators $x_1,...,x_m$ of the semigroup. This follows from the
relations. Hence if we can just take $m$ different maps of the
required nature, then they form a faithful representation. It is
certainly possible if $n \geq m$: take
$\as_1(j)=r_1,...,\as_m(j)=r_m, r_i \in N$. For different $r_i$ we
get different maps. \square

It is possible to find faithful representations of smaller
dimension. For example take  a subsets from theorem \ref{reps} of
size 3. Namely, let $m=3^d$ and our representation set consists of
the pairs $N=\{(k,\epsilon)| k=0,...,d-1, \epsilon= 0,1,2 \}$. Maps
define as follows: $\as_i (k,\epsilon)=(k,\epsilon_k(i))$, where
$\epsilon_k(i)$ is an $i$th coefficient in presentation of $i$ in
base $3$:
$i=\epsilon_0(i)+3\epsilon_1(i)+...+3^{d-1}\epsilon_{d-1}(i)$. We
have then a faithful representation on the set $N$ of size $3 \lceil
{\rm log}_3 m \rceil$.

It is not difficult to show that asymptotically this strategy gives
the best possible result, so asymptotically minimal size of faithful
representation is $3 {\rm log}_3 m$.

\subsection{Combinatorial description of maps corresponding to  maximal Hilbert series}

Here we give a combinatorial description of the maps $\sigma_i$,
$i=1,2$ satisfying condition (1) or (2) as they appear in the
theorem \ref{GrSig} above, that is those maps which define an
algebra with maximal Hilbert series.

 We also prove as a consequence  that if
Hilbert series is maximal then all maps $\sigma_i$, $i=\bar{1,m}$
coming from defining relations of arbitrary RIT algebra have
pairwise 2-isomorphic graphs.

Consider maps   $\sigma_i$ and $\sigma_k$. Suppose they obey
conditions described in theorem \ref{GrSig}. Let us define the set
$Y_0$ as a set of elements $j \in N$ where $\sigma_i$ and $\sigma_k$
coincide:

$$Y_0=\{j \in N | \sigma_i(j)=\sigma_k(j)\}.$$

\begin{lemma}\label{l1}
$Y_0$ is invariant under the action of both $\sigma_i$ and
$\sigma_k$
\end{lemma}

\proof Let $j$ be from $Y_0$. If in point $j$  $\sigma_i$ and
$\sigma_k$ obey condition (1) form theorem \ref{GrSig}, then
$\sigma_i(j) \in Y_0$ and $\sigma_k(j) \in Y_0$ due to the second
part of (1). Indeed, values of maps $\sigma_i$ and $\sigma_k$ on
their images should coincide, hence these images are again in $Y_0$.
Suppose in point $j$ condition (2) is fulfilled. Since $j \in Y_0$
we have $\sigma_i(j)=\sigma_k(j)=r$, but due to (2):
$r=\sigma_k(j)=\sigma_k\sigma_i(j)=\sigma_k(r)$ and
$r=\sigma_i(j)=\sigma_i\sigma_k(j)=\sigma_i(r)$. Which means that
actually (1) holds also for this point $j$. So for each point $j \in
Y_0$ condition (1) is satisfied. From this it easy follows that
$\sigma_i(Y_0)  \subset Y_0$ and $\sigma_k(Y_0)  \subset Y_0$.
\square

Let now consider an element $j \notin Y_0$ and its images
$j_1=\sigma_i(j)$ and $j_2=\sigma_k(j)$. Since images of $\sigma_i$
and $\sigma_k$ are different in $j$, condition (1) could not holds
in this point, thus we have condition (2) there. This gives us:
$j_2=\sigma_k(j)=\sigma_k \sigma_i (j)=\sigma_k(j_1)$ and
$j_1=\sigma_i(j)=\sigma_i \sigma_k (j)=\sigma_i(j_2)$. So we have
that $\sigma_k$ maps $j_1$ to $j_2$ and $\sigma_i$ maps $j_2$ to
$j_1$.

Using this information let us clarify how the element from outside
$Y_0$ could get to $Y_0$. Suppose $j \notin Y_0$ but
$j_1=\sigma_i(j) \in Y_0$. Then $\sigma_k(j_1)=\sigma_i(j_1)=j_2$.
This means that not only $j_2$ goes to $j_1$ under $\sigma_i$, but
also the other way around, $\sigma_i$ maps $j_1$ to $j_2$. On $Y_0$
condition (1) from the theorem \ref{GrSig} always holds and $j_1 \in
Y_0$ therefore for $j_2$ which is image of $j_1$ under $\sigma_i$ we
have $\sigma_i(j_2)=\sigma_k(j_2)$ (due to the second part of
condition (1) in point $j_1$). Since $j_1=\sigma_i(j_2) \in Y_0$,
$j_2$ is also in $Y_0$.

We have proved

\begin{lemma}\label{l2}
Let $j \notin Y_0$ but $\sigma_i(j) \in Y_0$. Then images of $j$
under $\sigma_i$ and $\sigma_k$ both are in $Y_0$ and  $\sigma_i$ as
well as $\sigma_k$ maps them to each other.

\end{lemma}

\begin{lemma}\label{l3}
If $j \notin Y_0$ but $j_1=\sigma_i(j) \in Y_0$, then there is no
such element from $N$, which has an image $j$ under $\sigma_i$ or
$\sigma_k$.
\end{lemma}

\proof Suppose there exists $m \in N$, such that say
$\sigma_i(m)=j$. Obviously $m \notin Y_0$, since $Y_0$ is invariant
and then $j$ should be in $Y_0$, but it is not. For points which are
not in $Y_0$ condition (1) could not holds, thus we have condition
(2) in $m$. This leads to the following contradiction:
$\sigma_k(m)=\sigma_k\sigma_i(m)=\sigma_k(j)=j_2,$
$\sigma_i(m)=\sigma_i\sigma_k(m)=\sigma_i(j_2)=j_1$ hence
$\sigma_i(m)=j$ and $\sigma_i(j)=j_1$, but $j_1$ can not be equal to
$j$ just because one is from $Y_0$ and another is not. \square

We are now in a position to define a bigger set $\tilde Y_0$:

$$\tilde Y_0 = Y_0 \cup \{j \in N | \sigma_i(j) \in Y_0  \},$$

which  satisfy the following nice property.

\begin{lemma}\label{split}
The set $N$ splits on two disjoint subsets which are invariant under
$\sigma_k$ and $\sigma_i$:
$$N=\tilde Y_0 \oplus P$$
\noindent where  $P=N \backslash \tilde Y_0$. Moreover the structure
of maps on $P$ is precisely as it was described in the theorem
\ref{reps}: $P$ is a disjoint union of subsets on which two points
are picked such that $\sigma_i$ maps entire  subset to one of them
and $\sigma_k$ to another.
\end{lemma}

\proof The fact that $\sigma_i$ and $\sigma_k$ preserve $Y_0$ was
proved above, invariance of $\tilde Y_0$ then follows from its
definition. Invariance of complement $P$ of $\tilde Y_0$ comes from
the statement of lemma \ref{l3}. \square

Let we define one more subset: $Z = Y_0 \backslash
\bigcup\limits_{j=i,k} \sigma_j(\tilde Y_0 \backslash Y_0)$.

Above lemmata allow us to give the following precise description of
maps $\sigma_i$, $\sigma_k$ corresponding to the maximal Hilbert
series.

\begin{theorem}\label{s1s2}
Algebra $R \in {\cal R}(2,n)$ has a maximal Hilbert series if and
only if maps $\sigma_i$, $\sigma_k$ coming from defining relations
has the following structure.

The set N is a disjoint union of invariant under both maps subsets
$P$ and $\tilde Y_0$:  $\,\, N=P \oplus \tilde Y_0$.

The action on $P$ is the following: $P$ is a disjoint union of
$P_i$, in each of them two points (not necessary different) are
fixed, such that $\sigma_i$ maps entire $P_i$ to one of them and
$\sigma_j$ to another.

The map on the other disjoint component $\tilde Y_0$ is the
following: there are three subsets $Z \subseteq Y_0 \subseteq \tilde
Y_0$. The set $Y_0 \backslash Z$ is a disjoint union of pairs $
\{j_1^{(l)},j_2^{(l)}\}_{l \in \Sigma}$, $\ (\Sigma$ is a finite set
of indexes),   such that $\forall j \in \tilde Y_0 \backslash Y_0
\,\, \exists l \in \Sigma: \sigma_i(j)=j_1^{(l)},
\sigma_k(j)=j_2^{(l)}$ and $\sigma_{i,k}(j_1^{(l)})=j_2^{(l)}$,
$\sigma_{i,k}(j_2^{(l)})=j_1^{(l)}$. Values of $\sigma_i$ and
$\sigma_k$ on $Z$ are in  $Y_0$ and coincide.

\end{theorem}

Let we say that two graphs are {\it 2-isomorphic} if their images
under gluing pairs of vertices of common cycles of length two are
isomorphic.

As a consequence of the above theorem we get the following

\begin{corollary}\label{cs1s2}
If an algebra $R \in {\cal R}(m,n)$ has a maximal Hilbert series,
then all graphs of maps $\sigma_i$ coming from the defining
relations of $R$ are pairwise 2-isomorphic.

\end{corollary}

\proof Note that the feature of condition on $\sigma_i$ formulated
in theorem \ref{GrSig} is that it is satisfied for an arbitrary set
$\{\sigma_i \}_{i=\bar {1,m}}$ if and only if it is satisfied for
any pair $\sigma_i$, $\sigma_k$, \,\, $i<k$. So we can apply theorem
\ref{s1s2} for any fixed pair $\sigma_i$, $\sigma_k$. It is clear
from the above description that graphs of maps $\sigma_i$ and
$\sigma_k$ are the same or isomorphic with the isomorphism defined
by permuting of some pairs $j_1,j_2$. Except one possible situation
when for some pair $j_1, j_2 \in Y_0 \backslash Z$ which consists of
images of $j \in \tilde Y_0 \backslash Y_0$ we have a point $l \in
Z$ such that $\sigma_i(l)= \sigma_k(l)=j_1$. This gives a
possibility for graphs of $\sigma_i$ and $\sigma_k$ to be
nonisomorphic. We exclude this possibility by gluing vertices of
common cycles of length two  ($j_1$ and $j_2$) in this two graphs,
so after such an operation graphs become isomorphic.  \square

\begin{corollary}\label{sl}
Combinatorial conditions from the theorem \ref{s1s2} on  $\sigma_i,
\sigma_k$ are equivalent to the following properties of algebra $R
\in {\cal R}(m,n) $:

(i) $R$ has a quadratic Gr\"obner basis;

(ii) $R$ has a lexicographically maximal in ${\cal R}(m,n) $ Hilbert
series;

(iii) $R$ is a PBW algebra (that is, have a series
$H_R=\frac{1}{(1-t)^{(m+n)}}$),

\noindent and imply the following properties of algebra:

(iv) $R$ is Koszul;

(v) $R$ is Auslander regular;

(vi) $R$ is Cohen-Macaulay.

\end{corollary}

\proof Equivalence to the condition (i) was proved in theorems
\ref{GrSig} and \ref{s1s2}. Equivalence of (i),(ii) and (iii) is a
direct consequence from the standard procedure of Hilbert series
computation for algebras presented by Gr\"obner basis. Implication
(i)$\Rightarrow$(iv) is known and could be found for example in
\cite{pio} or \cite{Gre}. Implications (i) $\Rightarrow$ (v) and (i)
$\Rightarrow$ (vi) could be found in \cite{wi}, \cite{LevAR},
\cite{Li}. \square

\section{Toward classification of Hilbert series}\label{tc}

Here we give a list of all RIT algebras of rank up to 4 with non
isomorphic colored graphs and the precise values of their Hilbert
series.

Denotation for a single graph (map) is the following: we write
$(i_1,i_2,...,i_k)$ for the map $\sigma:N \longrightarrow N: j
\mapsto i_j$.

We also denote by $ P_d(t)$ the series $\frac{1}{(1-t)^d}$, which is
a series of algebra $k[x_1,...,x_d]$ of commutative polynomials on
$d$ variables.

$\bf rk 1$ Commutative polynomials $k[x]$, with
$H_R(t)=P_1(t)=\frac{1}{1-t}$.

$\bf rk 2$

There are three possibilities: $R \in {\cal R}{(2,0)}$,   $R \in
{\cal R}{(1,1)}$, $R \in {\cal R}{(0,2)}$.

{\bf (2,0)}

$H_R(t)=P_2(t)=\frac{1}{(1-t)^2}$

{\bf (1,1)}

Graph: (1)

$H_R(t)=P_2(t)=\frac{1}{(1-t)^2}$

{\bf (0,2)}

$H_R(t)=P_2(t)=\frac{1}{(1-t)^2}$

$\bf rk 3$

There are four possibilities: $R \in {\cal R}{(3,0)}$, $R \in {\cal
R}{(2,1)}$, $R \in {\cal R}{(1,2)}$, $R \in {\cal R}{(0,3)}$.

{\bf (3,0)} and {\bf (0,3)}

$H_R(t)=P_3(t)=\frac{1}{(1-t)^3}$

{\bf (1,2)}

List of non isomorphic graphs: (1,2),(1,1),(2,1),(2,2)

All these four algebras have the series
$H_R(t)=P_3(t)=\frac{1}{(1-t)^3}$ due to the theorem \ref{reps}.

{\bf (2,1)}

2-colored graph: $\sigma_1=(1),\sigma_2=(1)$;
$H_R(t)=P_3(t)=\frac{1}{(1-t)^3}$ due to the theorem \ref{reps}.

$\bf rk 4$ There are five possibilities: $R \in {\cal R}{(4,0)}$, $R
\in {\cal R}{(3,1)}$, $R \in {\cal R}{(2,2)}$, $R \in {\cal
R}{(1,3)}$, $R \in {\cal R}{(0,4)}$.

{\bf (4,0)} and {\bf (0,4)}

$H_R(t)=P_4(t)=\frac{1}{(1-t)^4}$

{\bf (1,3)}

List of non isomorphic graphs:
(1,1,1),(1,1,2),(1,1,3),(1,2,3),(1,3,2),(2,1,1),(2,3,1)

All these algebras have the series $H_R(t)=P_4(t)=\frac{1}{(1-t)^4}$
due to the theorem \ref{reps}.

{\bf (2,2)}

List of non isomorphic 2-colored graphs:

 $\sigma_1=(1,2),\sigma_2=(1,2)$; $\sigma_1=(1,1),\sigma_2=(1,1)$;
 $\sigma_1=(1,1),\sigma_2=(2,2)$; $\sigma_1=(2,1),\sigma_2=(2,1)$

All algebras in this part of the list has maximal series:

$H_R(t)=P_4(t)=\frac{1}{(1-t)^4}$ due to the theorem \ref{reps}.

$\sigma_1=(1,2),\sigma_2=(1,1)$; $\sigma_1=(1,2),\sigma_2=(2,1)$;
 $\sigma_1=(1,1),\sigma_2=(2,1)$

 All algebras in this part of the list has series

 $H_R(t)=\frac{1+t+t^2}{1-3t+3t^2-t^3}$, which is not maximal since they does not
 obey conditions of the Corollary \ref{cs1s2}.

{\bf (3,1)}

3-colored graph: $\sigma_1=(1),\sigma_2=(1), \sigma_3=(1)$;
$H_R(t)=P_4(t)=\frac{1}{(1-t)^4}$ due to the theorem \ref{reps}.

\begin{corollary}\label{cl}
All RIT algebras of rank $\leq 4$ are PBW (with the series $P_d(t)$,
where $d$ is a rank) except of three cases, given by maps:

$R_1=\{\sigma_1=(1,2),\sigma_2=(1,1)\}$;

$R_2=\{\sigma_1=(1,2),\sigma_2=(2,1)\}$;

 $R_3=\{\sigma_1=(1,1),\sigma_2=(2,1)\}$.

\noindent These algebras have the series
$H_{R_i}(t)=\frac{1+t+t^2}{1-3t+3t^2-t^3}, \,\, i=\bar{1,3}$.

\end{corollary}

\section{Remark on the Generalized Yang--Baxter Equation for RIT algebras}

Once we have quadratic relations (\ref{fff}) of RIT algebra we could
define a linear map $r:V\otimes V\to V\otimes V$, where $V=\langle
x_1,\dots,x_m,y_1,\dots,y_n\rangle$ as follows:
\begin{align*}
r(x_i\otimes y_j)&=y_j\otimes x_i+\overline{y_{f(i,j)}\times y_j},
\\
r(x_i\otimes x_j)&=\overline{x_i\times x_j},
\\
r(y_i\otimes y_j)&=\overline{y_i\times y_j},
\end{align*}
where by $\overline{y_i\times y_j}$ we denote the product
$y_i\otimes y_j$ if $i<j$ and $y_j\otimes y_i$ if $i>j$. Analogously
$\overline{x_i\times x_j}$ stands for $x_i\otimes x_j$ if $i<j$ and
for $x_j\otimes x_i$ if $i>j$.

The map of this kind satisfies the Yang--Baxter equation if for the
action on $V\otimes V\otimes V$ induced by $r$ the following is true
$$
(r_{1,2}\otimes 1)(1\otimes r_{2,3})(r_{1,2}\otimes 1)= (1\otimes
r_{2,3})(r_{1,2}\otimes 1)(1\otimes r_{2,3}).
$$
Denote the operators on the left-hand and right-hand sides by
$R_{12}$ and $R_{23}$ respectively. Arguments of the same kind like
in the theorem \ref{GrSig} and \ref{s1s2}give us the following.

\begin{theorem} \label{th44}
If the relations (\ref{fff>}) of a RIT algebra $R$ form a Gr\"obner
basis, or equivalently if $\sigma_i$ satisfy the conditions of
Theorem \ref{GrSig}, then the corresponding operator $r:V\otimes V
\to V \otimes V$ defined as above from the relations, satisfies the
generalized Yang--Baxter equation $[R_{12},R_{23}]=0$.
\end{theorem}

Here in stead of the usual Yang--Baxter equation $R_{12}=R_{23}$ we
have got a generalized commutator version: $[R_{12},R_{23}]=0$.

In Theorem \ref{GrSig} we formulate combinatorial conditions on the
relations which mean that they form a quadratic Gr\"obner basis. In
Theorem \ref{th44} we state what kind of Yang--Baxter equation is
satisfied in this case. The latter has to do with the length of
chain of reductions necessary to ensure  that defining relations
form a Gr\"obner basis.

In more details connection to Yang--Baxter equations will be
discussed elsewhere, but we emphasize  here that combinatorial
description of maps in theorems \ref{GrSig} and \ref{s1s2} is also a
description of the case when defining relations may be associated to
a solutions of a kind of YBE.

\section{Acknowledgments}

This work was done during the stay of the second named author at the
Queen Mary, University of London as LMS Grace Young Fellow. In the
intermediate stage results were reported during the Special Semester
on Gr\"obner bases supported by RICAM (the Radon Institute for
Computational and Applied Mathematics, Austrian Academy of Sciences,
Linz) and RISC (Research Institute for Symbolyc Computation,
Johannes Kepler University, Linz, Austria) under the scientific
direction of Professor Bruno Buchberger. The second named author
also would like to thank the PITHAGORAS II program and the
Aristoteles University of Thessaloniki for hospitality and support
of this research.


\vspace{5 mm}

\begin{thebibliography}{99}

\itemsep=-2pt





\bibitem{An} D.~Anick, \it Generic algebras and CW complexes.
Algebraic topology and algebraic K-theory, \rm Proc. Conf.,
Princeton, NJ (USA), Ann. Math. Stud. 113(1987), 247-321.


\bibitem{wia}I.~Antoniou, N.~Iyudu, R.~Wisbauer, \it On Serre's
problem for  RIT Algebra, \rm Communications in Algebra
31(12)(2003), 6037-6050.

\bibitem{ia}I.~Antoniou, N.~Iyudu, \it Poincare--Hilbert
series, PI and Noetherianity of the enveloping of the RIT Lie
Algebra, \rm Communications in Algebra 29 (9)(2001), ~4183--4196.

\bibitem{a}I.~Antoniou,  \it Internal time and irreversibility of relativistic
dynamical systems, \rm PhD thesis, ULB, Brussels, 1988.

\bibitem{AS}V.~M.~Artin, W.~F.~Schelter, \it
Graded algebras of global dimension 3, \rm Adv. in Math. 66(1987),
171-216.

\bibitem{Bergm} G.~Bergman, \it The diamond lemma for ring theory.
\rm Adv. in Math. 29 (1978), no. 2, 178--218.

\bibitem{Buch} B.~Buchberger, {\it Ein Algorithmus zum Auffinden der Basiselemente des Restklassenringes
nach einem nulldimensionalen Polynomideal}, {1965}, Uni. Innsbruck.
(English translation: \it An Algorithm for Finding the Basis
Elements in the Residue Class Ring Modulo a Zero Dimensional
Polynomial Ideal, \rm Journal of Symbolic Computation, Special Issue
on Logic, Mathematics, and Computer Science: Interactions, 41(2006),
no.3-4, 475-511)

\bibitem{Gre} Ed.~Green, \it An introduction to noncommutative Gr\"obner bases, \rm LN
in Pure and Appl. Math., v.151, Dekker, New York, 1994

\bibitem{wi} N.~Iyudu, R.~Wisbauer, \it A counterexample to
Serre's conjecture in one class of quantum polynomial rings, \rm
Preprint, 2003.

\bibitem{prI} N.~Iyudu, {\it Classification of finite dimensional
representations of one noncommutative quadratic algebra},
Max-Planck-Institut f\"ur Mathematik Preprint series 2005, no. 20,
1-25.


\bibitem{LevG} V.~Levandovskyy, \it PBW Bases, Non--Degeneracy Conditions and
Applications, \rm Representation of algebras and related topics.
Proceedings of the ICRA X conference, v.45, AMS. Fields Institute
Communications, 2005, 229-246

\bibitem{LevAR} V.~Levandovskyy, \it Non--commutative Computer Algebra for polynomial algebras:
Gr\"obner bases, applications and implementation, \rm Doctoral
Thesis, Universit\"at Kaiserslautern, 2005,
http://kluedo.ub.uni-kl.de/volltexte/2005/1883/

\bibitem{Li} H.~Li, \it Noncommutative Gröbner bases and filtered-graded transfer.
\rm Lecture Notes in Mathematics, 1795. Springer-Verlag, Berlin,
2002, 197 pp


\bibitem{pio} D.~Piontkovski, \it Coherent algebras and noncommutative projective lines,
\rm Preprint, arXiv:math/0606279, 2006.

\bibitem{PP} A.~Polishchuk, L.~Positselski, \it Quadratic algebras,
\rm University Lecture Series, v 37, AMS, Provodence, RI, 2005.


\bibitem{TB} J.~Tate, M.~van den Bergh, \it Homological properties
of Sklyanin algebras, \rm Invent. Math., 124(1996), 619-647.

\bibitem{Ufn} V.~Ufnarovskij, \it Combinatorial and asymptotic
methods in algebra, Encycl. Math. Sci., v 57, Springer, Berlin 1995.




\end{thebibliography}
\end{document}